\documentstyle[11pt,twoside]{article}
\include{graphicx}
\title{The asymptotics of a generalised Struve function}
\author{\sc R. B.\ Paris \\
{\em Division of Computing and Mathematics}, \\
{\em Abertay University, Dundee DD1 1HG, UK}
}
\begin{document}
\def\f#1#2{\mbox{${\textstyle \frac{#1}{#2}}$}}
\def\dfrac#1#2{\displaystyle{\frac{#1}{#2}}}
\def\boldal{\mbox{\boldmath $\alpha$}}
{\newcommand{\Sgoth}{S\;\!\!\!\!\!/}
\newcommand{\bee}{\begin{equation}}
\newcommand{\ee}{\end{equation}}
\newcommand{\lam}{\lambda}
\newcommand{\ka}{\kappa}
\newcommand{\al}{\alpha}
\newcommand{\fr}{\frac{1}{2}}
\newcommand{\fs}{\f{1}{2}}
\newcommand{\g}{\Gamma}
\newcommand{\br}{\biggr}
\newcommand{\bl}{\biggl}
\newcommand{\ra}{\rightarrow}
\newcommand{\mbint}{\frac{1}{2\pi i}\int_{c-\infty i}^{c+\infty i}}
\newcommand{\mbcint}{\frac{1}{2\pi i}\int_C}
\newcommand{\mboint}{\frac{1}{2\pi i}\int_{-\infty i}^{\infty i}}
\newcommand{\gtwid}{\raisebox{-.8ex}{\mbox{$\stackrel{\textstyle >}{\sim}$}}}
\newcommand{\ltwid}{\raisebox{-.8ex}{\mbox{$\stackrel{\textstyle <}{\sim}$}}}
\renewcommand{\topfraction}{0.9}
\renewcommand{\bottomfraction}{0.9}
\renewcommand{\textfraction}{0.05}
\newcommand{\mcol}{\multicolumn}
\date{}
\maketitle
\pagestyle{myheadings}
\markboth{\hfill \sc R. B.\ Paris  \hfill}
{\hfill \sc  Asymptotics of the generalised Struve functions\hfill}
\begin{abstract}
A generalised Struve function has recently been introduced by Ali, Mondal and Nisar [J. Korean Math. Soc. {\bf 54} (2017) 575--598] as
\[(\fs z)^{\nu+1}\sum_{n=0}^\infty\frac{(\fs z)^{2n}}{\g(n+\f{3}{2}) \g(an+\nu+\f{3}{2})},\]
where $a$ is a positive integer. In this paper, we obtain the asymptotic expansions of this function for large complex $z$ when $a$ is a real parameter satisfying $a>-1$. Some numerical examples are presented to confirm the accuracy of the expansions. 
\vspace{0.4cm}

\vspace{0.3cm}

\end{abstract}

\vspace{0.3cm}

\noindent $\,$\hrulefill $\,$

\vspace{0.2cm}

\begin{center}
{\bf 1. \  Introduction}
\end{center}
\setcounter{section}{1}
\setcounter{equation}{0}
\renewcommand{\theequation}{\arabic{section}.\arabic{equation}}
The modified Struve function ${\bf L}_\nu(z)$ is a particular solution of the inhomogeneous modified Bessel equation
\[\frac{d^2y}{dz^2}+\frac{1}{z}\,\frac{dy}{dz}-\bl(1+\frac{\nu^2}{z^2}\br)y=\frac{(\fs z)^{\nu-1}}{\sqrt{\pi} \g(\nu+\fs)}\]
which possesses the series expansion \cite[p.~288]{DLMF}
\bee\label{e11}
{\bf L}_\nu(z)=(\fs z)^{\nu+1} \sum_{n=0}^\infty \frac{(\fs z)^{2n}}{\g(n+\f{3}{2}) \g(n+\nu+\f{3}{2})}
\ee
valid for all finite $z$. The ordinary Struve function ${\bf H}_\nu(z)$ is given by the alternating version of (\ref{e11}) corresponding to $z$ situated on the imaginary axis, and
\bee\label{e12}
{\bf H}_\nu(z)=\mp i e^{\mp\pi i\nu/2} {\bf L}_\nu(\pm iz).
\ee

A generalisation of the series (\ref{e11}) and its alternating version has been considered by Yagmur and Orban  \cite{YO} by replacing the $\f{3}{2}$ in the second gamma function by an arbitrary
complex parameter. These authors determined sufficient conditions for it to be univalent and obtained various convexity properties of the functions.
In a recent paper by Ali, Mondal and Nisar \cite{AMN}, the above series was generalised further by the introduction of 
a multiple argument in the second gamma function, namely the function 
\[\sum_{n=0}^\infty\frac{(\pm 1)^n (\fs z)^{2n+\nu+1}}{\g(n+\f{3}{2}) \g(an+\nu+\f{3}{2})},\]
where $a$ denotes a positive integer. These authors showed that this function is a solution of an $(a+1)$th-order differential equation. They also investigated its monotonicity and log-convexity properties and established Tur\'an-type inequalities. Upper bounds satisfied by this function in the case $a=2$ and for the non-alternating series were derived.

In this paper we consider the generalised Struve function defined in \cite{AMN} by
\bee\label{e13}
{\bf L}_\nu(z;a)=(\fs z)^{\nu+1} \sum_{n=0}^\infty \frac{(\fs z)^{2n}}{\g(n+\f{3}{2}) \g(an+\nu+\f{3}{2})}
\ee
where $a$ is a real parameter satisfying $a>-1$ for convergence of the series and we restrict $\nu$ to be real.
It is then evident that
\[{\bf L}_\nu(ze^{\pm\pi i};a)=e^{\pm\pi i(\nu+1)} {\bf L}_\nu(z;a),\qquad {\bf L}_\nu({\overline z};a)=\overline{{\bf L}_\nu(z;a)},\]
where the bar denotes the complex conjugate,
so that we may confine our attention to the sector $0\leq \arg\,z\leq\fs\pi$. When $\arg\,z=\fs\pi$ in (\ref{e13}), we generate the generalised alternating version given by ${\bf H}_\nu(|z|;a)$. 
We determine the asymptotic expansion of ${\bf L}_\nu(z;a)$ for large complex $z$ and finite values of $\nu$ and $a$. The series in (\ref{e13}) is a particular case of a generalised Wright function; see (\ref{e20}) below.
Accordingly, we employ the well-established asymptotic theory of the Wright function to determine the large-$z$ expansion of ${\bf L}_\nu(z;a)$, a summary of which is presented in Section 2. It will be found that the analysis
of ${\bf L}_\nu(z;a)$ separates into two distinct cases according as $a>0$ and $-1<a<0$.
\vspace{0.6cm}

\begin{center}
{\bf 2. \ Standard asymptotic theory of the generalised Wright function}
\end{center}
\setcounter{section}{2}
\setcounter{equation}{0}
\renewcommand{\theequation}{\arabic{section}.\arabic{equation}}
The generalised Wright function is defined by the series
\bee\label{e20}
{}_p\Psi_q(z)\equiv{}_p\Psi_q\bl(\!\!\begin{array}{c}(\alpha_1,a_1), \ldots ,(\alpha_p,a_p)\\(\beta_1, b_1), \ldots ,(\beta_q,b_q)\end{array}\!\!\bl|z\!\br)=\sum_{n=0}^\infty g(n)\,\frac{z^n}{n!},
\ee
\bee\label{e20a}
g(n)=\frac{\prod_{r=1}^p\Gamma(\alpha_rn+a_r)}{\prod_{r=1}^q\Gamma(\beta_rn+b_r)},
\ee
where $p$ and $q$ are nonnegative integers, the parameters $\alpha_r$  and 
$\beta_r$ are real and positive and $a_r$ and $b_r$ are
arbitrary complex numbers. We also assume that the $\alpha_r$ and $a_r$ are subject to 
the restriction
\bee\label{e20ab}
\alpha_rn+a_r\neq 0, -1, -2, \ldots \qquad (n=0, 1, 2, \ldots\ ;\, 1\leq r \leq p)
\ee
so that no gamma function in the numerator in (\ref{e20}) is singular.

We introduce the parameters associated\footnote{Empty sums and products are to be interpreted as zero and unity, respectively.} with $g(n)$ given by
\[\kappa=1+\sum_{r=1}^q\beta_r-\sum_{r=1}^p\alpha_r, \qquad 
h=\prod_{r=1}^p\alpha_r^{\alpha_r}\prod_{r=1}^q\beta_r^{-\beta_r},\]
\bee\label{e21}
\vartheta=\sum_{r=1}^pa_r-\sum_{r=1}^qb_r+\f{1}{2}(q-p),\qquad \vartheta'=1-\vartheta.
\ee
If it is supposed that $\alpha_r$ and $\beta_r$ are such that $\kappa>0$ then ${}_p\Psi_q(z)$ 
is uniformly and absolutely convergent for all finite $z$. If $\kappa=0$, the sum in (\ref{e11})
has a finite radius of convergence equal to $h^{-1}$, whereas for $\kappa<0$ the sum is divergent 
for all nonzero values of $z$. The parameter $\kappa$ will be found to play a critical role 
in the asymptotic theory of ${}_p\Psi_q(z)$ by determining the sectors in the $z$-plane 
in which its behaviour is either exponentially large, algebraic or exponentially small 
in character as $|z|\ra\infty$.

The determination of the asymptotic expansion of ${}_p\Psi_q(z)$ for $|z|\ra\infty$ and finite 
values of the parameters has a long history.
Detailed investigations were carried out by Wright \cite{W1, W2} and by
Braaksma \cite{Br} for a more general class of integral functions than (\ref{e11}). We present below a summary of the main expansion theorems related to the asymptotics of ${}_p\Psi_q(z)$ for large $|z|$; for a recent presentation, see \cite{P17}. In order to do this
we first introduce the exponential expansion $E_{p,q}(z)$ and the 
algebraic expansion $H_{p,q}(z)$ associated with ${}_p\Psi_q(z)$.

The exponential expansion $E_{p,q}(z)$ is given by the formal asymptotic sum
\bee\label{e22c}
E_{p,q}(z):=Z^\vartheta e^Z\sum_{j=0}^\infty A_jZ^{-j}, \qquad Z=\kappa (hz)^{1/\kappa},
\ee
where the coefficients $A_j$ are those appearing in the inverse factorial expansion of $g(s)/s!$ given by  
\bee\label{e22a}
\frac{g(s)}{\g(1+s)}=\kappa (h\kappa^\kappa)^s\bl\{\sum_{j=0}^{M-1}\frac{A_j}{\Gamma(\kappa s+\vartheta'+j)}
+\frac{\rho_M(s)}{\Gamma(\kappa s+\vartheta'+M)}\br\}.
\ee
Here $g(s)$ is defined in (\ref{e20a}) with $n$ replaced by $s$, $M$ is a positive integer and $\rho_M(s)=O(1)$ for $|s|\ra\infty$ in $|\arg\,s|<\pi$.
The leading coefficient $A_0$ is specified by
\bee\label{e22b}
A_0=(2\pi)^{\frac{1}{2}(p-q)}\kappa^{-\frac{1}{2}-\vartheta}\prod_{r=1}^p
\alpha_r^{a_r-\frac{1}{2}}\prod_{r=1}^q\beta_r^{\frac{1}{2}-b_r}.
\ee
The coefficients $A_j$ are independent of $s$ and depend only on the parameters $p$, $q$, $\alpha_r$, 
$\beta_r$, $a_r$ and $b_r$. An algorithm for their evaluation is described in \cite{P17, P17a}.

The algebraic expansion $H_{p,q}(z)$ follows from the Mellin-Barnes integral representation \cite[\S 2.4]{PK}
\bee\label{e24aa}
{}_p\Psi_q(z)=\frac{1}{2\pi i}\int_{-\infty i}^{\infty i} \Gamma(s)g(-s)(ze^{\mp\pi i})^{-s}ds,\qquad |\arg(-z)|<\pi(1-\fs\kappa),
\ee
where the path of integration is indented near $s=0$ to separate\footnote{This is always 
possible when the condition (\ref{e20ab}) is satisfied.} the poles of $\g(s)$ at $s=-k$ from those of 
$g(-s)$ situated at 
\bee\label{e24a}
s=(a_r+k)/\alpha_r, \qquad k=0, 1, 2, \dots\, \ (1\leq r\leq p).
\ee
In general there will be $p$ such sequences of simple poles though, depending on the values 
of $\alpha_r$ and $a_r$, some of these poles could be multiple poles or even ordinary 
points if any of the $\Gamma(\beta_rs+b_r)$ are singular there. Displacement of the contour to the 
right over the poles of $g(-s)$ then yields the algebraic expansion of 
${}_p\Psi_q(z)$ valid in the sector in (\ref{e24aa}). 

If it is assumed that the parameters are such that 
the poles in (\ref{e24a}) are all simple we obtain the algebraic expansion given by 
$H_{p,q}(z)$, where
\bee\label{e25}
H_{p,q}(z):=\sum_{m=1}^p\alpha_m^{-1}z^{-a_m/\alpha_m}S_{p,q}(z;m)
\ee
and $S_{p,q}(z;m)$ denotes the formal asymptotic sum
\bee\label{e25a}
S_{p,q}(z;m):=\sum_{k=0}^\infty \frac{(-)^k}{k!}\Gamma\left(\frac{k+a_m}{\alpha_m}\right)\,
\frac{\prod_{r=1}^{'\,p}\Gamma(a_r-\alpha_r(k+a_m)/\alpha_m)}
{\prod_{r=1}^q\Gamma(b_r-\beta_r(k+a_m)/\alpha_m)} z^{-k/\alpha_m},
\ee
with the prime indicating the omission of the term corresponding to $r=m$ in the product. 
This expression in (\ref{e25}) consists of (at most) $p$ expansions each with the leading behaviour 
$z^{-a_m/\alpha_m}$ ($1\leq m\leq p$).
When the parameters $\alpha_r$ and $a_r$ are such that some of the poles 
are of higher order, the expansion (\ref{e25a}) is invalid and the residues must 
then be evaluated according to the multiplicity of the poles concerned; this will lead to terms involving $\log\,z$ in the algebraic expansion.

The expansion theorems for ${}_p\Psi_q(z)$ are as follows. Throughout we let $\epsilon$ denote an arbitrarily 
small positive quantity.
\newtheorem{theorem}{Theorem}
\begin{theorem}$\!\!\!.$
When $0<\kappa\leq 2$, then 
\bee\label{e24}
{}_p\Psi_q(z)\sim\left\{\begin{array}{lll}
E_{p,q}(z)+H_{p,q}(ze^{\mp\pi i}) & 
\mbox{in} & |\arg\,z|\leq \min\{\pi-\epsilon,\pi\kappa-\epsilon\} \\
\\H_{p,q}(ze^{\mp\pi i}) & 
\mbox{in} & \pi\kappa+\epsilon\leq |\arg\,z|\leq\pi\ \ (0<\kappa<1)\\
\\E_{p,q}(z)+E_{p,q}(ze^{\mp2\pi i})+H_{p,q}(ze^{\mp\pi i}) & 
\mbox{in} & |\arg\,z|\leq\pi\ \ (1<\kappa\leq 2)\end{array} \right.
\ee
as $|z|\ra\infty$. The upper or lower signs 
are chosen according as $\arg\,z>0$ or $\arg\,z<0$, respectively. 
\end{theorem}
\begin{theorem}$\!\!\!.$
When $\kappa>2$ we have
\bee\label{e23}
{}_p\Psi_q(z)\sim \sum_{n=-N}^N E_{p,q}(ze^{2\pi in})+H_{p,q}(ze^{\mp\pi iz})
\ee
as $|z|\to\infty$ in the sector $|\arg\,z|\leq\pi$. The integer $N$ is chosen such that it is the smallest integer satisfying $2N+1>\fs\kappa$ and the upper or lower is chosen according as $\arg\,z>0$ or $\arg\,z<0$, respectively.
\end{theorem}
In this case the asymptotic behaviour of ${}_p\Psi_q(z)$ is exponentially large for all values of $\arg\,z$ and, consequently, the algebraic expansion may be neglected. The sums $E_{p,q}(ze^{2\pi in})$ are exponentially large (or oscillatory) as $|z|\to\infty$ for values of $\arg\,z$ satisfying $|\arg\,z+2\pi n|\leq\fs\pi\kappa$.

The division of the $z$-plane into regions where ${}_p\Psi_q(z)$ possesses exponentially large or algebraic behaviour for large $|z|$ is illustrated in Fig.~1.
When $0<\kappa<2$, the exponential expansion $E_{p,q}(z)$ is still present in the sectors $\fs\pi\kappa<|\arg\,z|<\min\{\pi,\pi\kappa\}$, where it is subdominant. The rays $\arg\,z=\pm\pi\kappa$ ($0<\kappa<1$), where $E_{p,q}(z)$ is {\it maximally} subdominant with respect to $H_{p,q}(ze^{\mp\pi i})$, are called Stokes lines.\footnote{The positive real axis $\arg\,z=0$ is also a Stokes line where the algebraic expansion is maximally subdominant.} As these rays are crossed (in the sense of increasing $|\arg\,z|$) the exponential expansion switches off according to Berry's now familiar error-function smoothing law \cite{B}; see \cite{P10} for details. The rays $\arg\,z=\pm\fs\pi\kappa$, where $E_{p,q}(z)$ is oscillatory and comparable to $H_{p,q}(ze^{\mp\pi i})$, are called anti-Stokes lines.
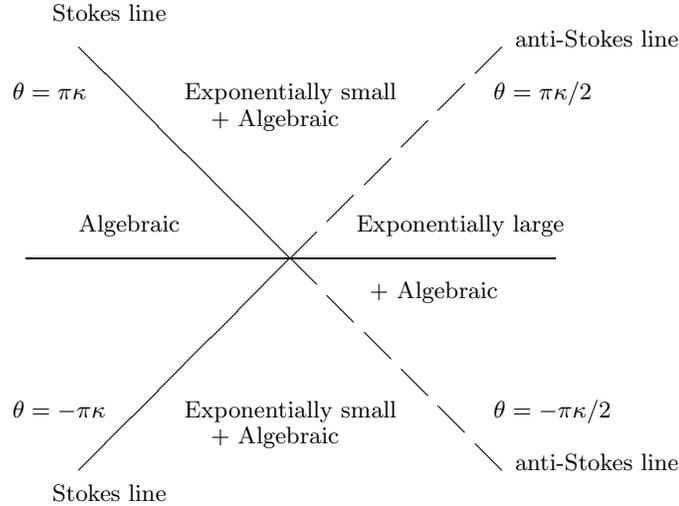
\begin{figure}[t]
\centering
\begin{picture}(200,200)(0,0)
\put(0,100){\line(1,0){200}}
\put(100,100){\line(-1,1){80}}
\put(100,100){\line(-1,-1){80}}
\multiput(100,100)(14,14){6}{\line(1,1){10}}
\multiput(100,100)(14,-14){6}{\line(1,-1){10}}
\put(125,110){\footnotesize{Exponentially large}}
\put(130,85){\footnotesize{$+$ Algebraic}}
\put(60,160){\footnotesize{Exponentially small}}
\put(70,150){\footnotesize{$+$ Algebraic}}
\put(60,40){\footnotesize{Exponentially small}}
\put(70,30){\footnotesize{$+$ Algebraic}}
\put(20,110){\footnotesize{Algebraic}}
\put(177,160){\footnotesize{$\theta=\pi\kappa/2$}}
\put(185,180){\footnotesize{anti-Stokes line}}
\put(177,40){\footnotesize{$\theta=-\pi\kappa/2$}}
\put(185,20){\footnotesize{anti-Stokes line}}
\put(-5,160){\footnotesize{$\theta=\pi\kappa$}}
\put(10,190){\footnotesize{Stokes line}}
\put(-5,40){\footnotesize{$\theta=-\pi\kappa$}}
\put(10,8){\footnotesize{Stokes line}}
\end{picture}
\caption{\small{The exponentially large and algebraic sectors associated with ${}_p\Psi_q(z)$ in the complex $z$-plane with $\theta=\arg\,z$ when $0<\kappa<1$. The Stokes and anti-Stokes lines are indicated.}}
\end{figure}

We omit the expansion {\it on} the Stokes lines $\arg\,z=\pm\pi\kappa$ in Theorem 1; the details in the case $p=1$, $q\geq0$ are discussed in \cite{P14}; see also \cite{P17a} for the case of the generalised Bessel function.
Since $E_{p,q}(z)$ is exponentially small in $\fs\pi\kappa<|\arg\,z|\leq\pi$, then in the sense of Poincar\'e, the expansion $E_{p,q}(z)$ can be neglected in these sectors.
Similarly, $E_{p,q}(ze^{-2\pi i})$ is exponentially small compared to $E_{p,q}(z)$ in $0\leq\arg\,z<\pi$ and so can be neglected when $1<\kappa<2$. However, in the vicinity of $\arg\,z=\pi$, these last two expansions are of comparable magnitude and, for real parameters, they combine to generate a real result on this ray. A similar remark applies to $E_{p,q}(ze^{2\pi i})$ in $-\pi<\arg\,z\leq 0$. 
\vspace{0.6cm}

\begin{center}
{\bf 3. \ The asymptotic expansion of ${\bf L}_\nu(z;a)$ when $a>0$ }
\end{center}
\setcounter{section}{3}
\setcounter{equation}{0}
\renewcommand{\theequation}{\arabic{section}.\arabic{equation}}
The function ${\bf L}_\nu(z;a)$ can be written as
\bee\label{e31}
(\fs z)^{-\nu-1}\,{\bf L}_\nu(z;a)={}_1\Psi_2\bl(\begin{array}{c} (1,1)\\(1,\f{3}{2}), (a,\nu+\f{3}{2})\end{array}\bl |\,\zeta\br)\equiv {}_1\Psi_2(\zeta),\qquad \zeta=\f{1}{4}z^2,
\ee
where the Wright function ${}_1\Psi_2(\zeta)$ is defined in (\ref{e20}).
From (\ref{e21}) and (\ref{e22b}) we have the parameters associated with the right-hand side of (\ref{e31})
\[\kappa=1+a,\qquad h=a^{-a},\qquad \vartheta=-\nu-\f{3}{2},\qquad A_0=\frac{1}{\sqrt{2\pi}}\,(\kappa/a)^{\nu+1}.\]
We note that $\kappa>1$ when $a>0$. 

The algebraic and exponential expansions associated with ${}_1\Psi_2(\zeta)$ are, from (\ref{e22c}) and (\ref{e25}),
\[H_{1,2}(\zeta)=\frac{1}{\pi} \sum_{k=0}^\infty \frac{\g(k+\fs)\,\zeta^{-k-1}}{\g(\nu+\f{3}{2}-a(1+k))}\]
and 
\[E_{1,2}(\zeta)=A_0Z^\vartheta e^Z \sum_{j=0}^\infty c_j Z^{-j},\qquad  Z=\kappa(h\zeta)^{1/\kappa},\]
where $c_0=1$. The coefficients $c_j\equiv c_j(a,\nu)$ ($j\geq 1$) can be determined by the algorithm described in \cite[Appendix]{P17}; see also \cite{P17a}, \cite[p.~46]{PK}. It is found that
\[c_1(a,\nu)=-\frac{1}{24a}\{11+24\nu+12\nu^2-a(25+24\nu)+11a^2\},\]
\[c_2(a,\nu)=\frac{1}{1152a^2}\{265+1056\nu+1416\nu^2+768\nu^3+144\nu^4-2a(791+2040\nu+1596\nu^2+384\nu^3)\]\[+3a^2(905+1360\nu+472\nu^2)-2a^3(791+528\nu)+265a^4\},\]
\[c_3(a,\nu)=-\frac{1}{414720a^3}\{(48703+286200\nu+617940\nu^2+636480\nu^3+334800\nu^4+86400\nu^5+8640\nu^6)\]
\[-3a(189797 + 791400\nu + 1179240\nu^2 + 797760\nu^3 + 248400\nu^4 + 
 28800\nu^5)\]
 \[+6 a^2 (355459 + 1019700\nu + 996570\nu^2 + 398880\nu^3 + 55800\nu^4)\]
 \[-a^3 (3254507 + 6118200\nu + 3537720\nu^2 + 636480\nu^3)\]
 \bee\label{e31a}
 +6 a^4 (355459 + 395700\nu + 102990\nu^2)-3 a^5 (189797 + 95400\nu)+48703 a^6\}.
 \ee
 The complexity of the higher coefficients prevents their presentation in general form. However, in specific cases, where the parameters have numerical values, it is possible to generate many coefficients; see the example in Section 5.
 
Then from Theorems 1 and 2 we obtain 
\begin{theorem}$\!\!\!.$\ When $a>0$ and $\zeta=z^2/4$, we have the expansion of the generalised Struve function 
\bee\label{e32}
(\fs z)^{-\nu-1} {\bf L}_\nu(z;a)\sim\left\{\begin{array}{ll} E_{1,2}(\zeta)+E_{1,2}(\zeta e^{\mp2\pi i})+H_{1,2}(\zeta e^{\mp\pi i}) & (1<\kappa\leq 2)\\
\\
\displaystyle{\sum_{n=-N}^N E_{1,2}(\zeta e^{2\pi in})+H_{1,2}(\zeta e^{\mp\pi i})} & (\kappa>2)\end{array}\right.
\ee
as $|z|\to\infty$ in the sector $|\arg\,z|\leq\fs\pi$. The integer $N$ is chosen such that it is the smallest integer satisfying $2N+1>\fs\kappa$ and the upper or lower signs are chosen according as $\arg\,z>0$ or $\arg\,z<0$, respectively. 
\end{theorem}
The expansion $E_{1,2}(\zeta e^{-2\pi i})$ is exponentially smaller than $E_{1,2}(\zeta)$ in $0\leq\arg\,\zeta<\pi$ ($0\leq\arg\,z<\fs\pi$) and can be neglected,
depending on the accuracy required, but becomes comparable to $E_{1,2}(\zeta)$ in the vicinity of $\arg\,\zeta=\pi$. A similar remark applies to $E_{1,2}(\zeta e^{2\pi i})$ in the lower half-plane.
In fact, the expansions $E_{1,2}(\zeta e^{\mp2\pi i})$ switch off as $|\arg\,\zeta|$ decreases across the Stokes lines for these
functions; see \cite[\S 3.4]{P17} for a numerical example. Since the exponential factors associated with $E_{1,2}(\zeta)$ and $E_{1,2}(\zeta e^{\mp2\pi i})$ are $\exp\,[|Z|e^{i\phi/\kappa}]$ and $\exp\, [|Z|e^{i(\phi\mp2\pi)/\kappa}]$, where $\phi=\arg\,\zeta$, the greatest difference between the real parts of these factors occurs when $\sin (\phi/\kappa)=\sin ((\phi\mp2\pi)/\kappa)$; that is, on the Stokes lines $\phi=\pm\pi(1-\fs\kappa)$. Consequently, when $1<\kappa<2$, the expansions $E_{1,2}(\zeta e^{\mp2\pi i})$ are not present as $|z|\to\infty$ in the sector $|\arg\,z|<\fs\pi(1-\fs\kappa)$.

When $z>0$ ($\arg\,\zeta=0$), the algebraic expansion in (\ref{e32}) is maximally subdominant as $z\to+\infty$ and undergoes a Stokes phenomenon. Neglecting this subdominant contribution, we therefore have
\bee\label{e33a}
(\fs z)^{-\nu-1} {\bf L}_\nu(z;a)\sim E_{1,2}(\zeta)\qquad (1<\kappa\leq 2;\ z\to+\infty).
\ee
When $\arg\,z=\fs\pi$ ($\arg\,\zeta=\pi$) we have from (\ref{e32}) with $\zeta=xe^{\pi i}$ the expansion
\[(\fs i|z|)^{-\nu-1} {\bf L}_\nu(i|z|;a)\sim E_{1,2}(xe^{\pi i})+E_{1,2}(xe^{-\pi i})+H_{1,2}(x)\]
\bee\label{e33}
={\tilde E}_1(X)+H_{1,2}(x)\qquad (1<\kappa\leq 2)
\ee
as $|z|\to+\infty$, where 
\[{\tilde E}_n(X):=2A_0X^\vartheta e^{X \cos\,(2n-1)\pi/\kappa} \sum_{j=0}^\infty c_jX^{-j} \cos\,\bl[X\sin\,\frac{(2n-1)\pi}{\kappa} +\frac{\pi}{\kappa}(\vartheta-j)\br]\]
for positive integer $n$
and $$x=\f{1}{4} |z|^2,\qquad X=\kappa(hx)^{1/\kappa}.$$ The algebraic expansion $H_{1,2}(x)$ is dominant as $|z|\to\infty$.

When $\kappa>2$, the expansion of $(\fs z)^{-\nu-1} {\bf L}_\nu(z;a)$ as $z\to+\infty$ is given by the second expression in (\ref{e32}), where the maximally subdominant algebraic expansion may be neglected. When $z=i|z|$, we have
\bee\label{e34}
(\fs i|z|)^{-\nu-1} {\bf L}_\nu(i|z|;a)\sim\sum_{n=1}^N {\tilde E}_n(X)+H_{1,2}(x)
\ee
as $|z|\to\infty$, 
where $N$ is as specified above and we have omitted the exponentially small contribution represented by $E_{1,2}(xe^{(2N+1)\pi i})$. The growth of ${\bf L}_\nu(i|z|;a)$ in this case is exponentially large as $|z|\to\infty$.

\vspace{0.6cm}

\begin{center}
{\bf 4. \ The asymptotic expansion of ${\bf L}_\nu(z;a)$ when $-1<a<0$ }
\end{center}
\setcounter{section}{4}
\setcounter{equation}{0}
\renewcommand{\theequation}{\arabic{section}.\arabic{equation}}
When $-1<a<0$ we write $a=-\sigma$, where $0<\sigma<1$. Then use of the reflection formula for the gamma function shows that
\begin{eqnarray}
(\fs z)^{-\nu-1} {\bf L}_\nu(z;-\sigma)&=&\frac{1}{\pi}\sum_{n=0}^\infty \frac{\g(\sigma n-\nu-\fs) \,(\fs z)^{2n}}{\g(n+\f{3}{2})}\,\sin\,\pi(-\sigma n+\nu+\f{3}{2})\nonumber\\
&=&\frac{i}{2\pi}\bl\{e^{\pi i\vartheta} {}_2\Psi_1(\zeta e^{\pi i\sigma})-e^{-\pi i\vartheta} {}_2\Psi_1(\zeta e^{-\pi i\sigma})\br\},\label{e41}
\end{eqnarray}
where $\zeta=z^2/4$ and $\vartheta=-\nu-\f{3}{2}$ as in Section 3 and
\bee\label{e41b}
{}_2\Psi_1(\zeta)\equiv{}_2\Psi_1\bl(\begin{array}{c} (1,1), (\sigma,-\nu-\fs)\\(1,\f{3}{2})\end{array}\bl |\,\zeta\br)=\sum_{n=0}^\infty \frac{\g(\sigma n\!-\!\nu\!-\!\fs) \,\zeta^n}{\g(n+\f{3}{2})}
\ee
provided $\g(\sigma n-\nu-\fs)$ is regular for $n=0, 1, 2, \ldots\ $.

The parameters associated with the Wright function in (\ref{e41b}) are
\bee\label{e41a}
\kappa=1-\sigma,\qquad h=\sigma^\sigma,\qquad A_0=\sqrt{2\pi}\,(\kappa/\sigma)^{\nu+1}.
\ee
Then, since $0<\kappa<1$, we have from Theorem 1
\bee\label{e42}
{}_2\Psi_1(\zeta)\sim \left\{\begin{array}{ll} E_{2,1}(\zeta)+H_{2,1}(\zeta e^{\mp\pi i}) & \mbox{in}\ |\arg\,\zeta|\leq\pi\kappa-\epsilon\\
\\
H_{2,1}(\zeta e^{\mp\pi i}) & \mbox{in}\ \pi\kappa+\epsilon\leq |\arg\,\zeta|\leq\pi\end{array}\right.
\ee
as $\zeta\to\infty$, where  from (\ref{e22c}) the exponential expansion is
\bee\label{e43}
E_{2,1}(\zeta)=A_0 Z^\vartheta e^Z \sum_{j=0}^\infty c_j Z^{-j},\qquad Z=\ka(h\zeta)^{1/\kappa},
\ee
with the coefficients $c_j=c_j(-\sigma,\nu)$ given in (\ref{e31a}); see the appendix. From (\ref{e25}), the algebraic expansion consists of two asymptotic series, viz.
\bee\label{e44}
H_{2,1}(\zeta)=\sum_{k=0}^\infty (-)^kG_k \zeta^{-k-1}+\frac{1}{\sigma} \sum_{k=0}^\infty (-)^kG'_k \zeta^{-k_s},\quad k_s:=\frac{1}{\sigma}(k-\nu-\fs),
\ee
where, provided $k_s$ is not equal to an integer\footnote{When $k=1, 2, \ldots\ $ there are double poles in the integrand of (\ref{e24aa}); the values $k_s=0, -1, -2, \ldots\ $ are disallowed by (\ref{e24a}).},
\[G_k=\frac{\g(-\nu\!-\!\fs\!-\!\sigma(k\!+\!1))}{\g(\fs-k)},\qquad G_k'=\frac{\g(k_s)\g(1-k_s)}{k!\,\g(\f{3}{2}-k_s)}.\]
The expansion of ${}_2\Psi_1(\zeta)$ on the Stokes lines $\arg\,\zeta=\pm\pi\kappa$, where the exponential expansion $E_{2,1}(\zeta)$ is maximally subdominant and is in the process of switching off (as $|\arg\,\zeta|$ increases), is omitted here. This has been considered for the more general function ${}_1\Psi_q(z)$ for integer $q\geq0$ in \cite[\S 5]{P14} and for the generalised Bessel function ${}_0\Psi_1(z)$ in \cite{P17a}. 

The expansion of ${\bf L}_\nu(z;-\sigma)$ can be constructed from knowledge of the expansion of the associated function ${}_2\Psi_1(\zeta)$ in (\ref{e42}). To keep the presentation as clear as possible, we restrict our attention here to the most commonly occurring situation of real $\zeta$ (that is, $z>0$ or $z=i|z|$) in (\ref{e41}).
\bigskip

\noindent 4.1 \ {\it The expansion for $z\to+\infty$}
\vspace{0.2cm}

\noindent 
When $z>0$ ($\zeta>0$) we write $\zeta=x$. From (\ref{e41}) and (\ref{e44}), the algebraic component of the expansion of $(\fs z)^{-\nu-1} {\bf L}_\nu(z;-\sigma)$ is
\[{\hat H}_{2,1}(x)\equiv\frac{i}{2\pi}\bl\{e^{\pi i\vartheta} H_{2,1}(xe^{\pi i\sigma}\,\cdot e^{-\pi i})-e^{-\pi i\vartheta} H_{2,1}(xe^{-\pi i\sigma}\,\cdot e^{\pi i})\br\}\]
\[=\frac{1}{\pi}\sum_{k=0}^\infty G_k x^{-k-1} \sin \pi(\vartheta-\sigma(k+1))-\frac{1}{\pi\sigma}\sum_{k=0}^\infty(-)^k G_k' x^{-k_s} \sin \pi(\vartheta+\kappa k_s)\]
\bee\label{e45}
=\frac{1}{\pi}\sum_{k=0}^\infty \frac{(-x)^{-k-1} \g(k+\fs)}{\g(\nu\!+\!\f{3}{2}\!+\!\sigma(k+1))}+\frac{1}{\sigma}\sum_{k=0}^\infty \frac{x^{-k_s}}{k! \g(\f{3}{2}-k_s)}.
\ee
In (\ref{e41b}) it was necessary to assume that $\g(\sigma n-\nu-\fs)$ is regular for $n=0, 1, 2, \ldots\ $. If this assumption is false, the second expansion appearing in (\ref{e45}) is still valid; compare \cite[\S 4]{W5}.
The exponential component is from (\ref{e43}) given by
\[{\hat E}_{2,1}(x)\equiv \frac{i}{2\pi}\bl\{e^{\pi i\vartheta} E_{2,1}(xe^{\pi i\sigma})-e^{-\pi i\vartheta} E_{2,1}(xe^{-\pi i\sigma})\br\}\]
\bee\label{e46}
=\frac{A_0}{\pi} X^\vartheta e^{X\cos\,(\pi\sigma/\kappa)} \sum_{j=0}^\infty (-)^{j-1}c_j X^{-j} \sin \bl[X\sin \frac{\pi\sigma}{\kappa}+\frac{\pi}{\kappa}(\vartheta-j)\br],
\ee
where $X=\kappa(hx)^{1/\kappa}$.
\begin{figure}[t]
\centering
\begin{picture}(300,200)(0,0)
\put(0,100){\line(1,0){120}}
\put(180,100){\line(1,0){120}}
\put(60,100){\line(-1,1){60}}
\put(60,100){\line(-1,-1){60}}
\put(60,100){\line(3,1){70}}
\put(60,100){\line(3,-1){70}}
\put(240,100){\line(-1,1){60}}
\put(240,100){\line(-1,-1){60}}
\multiput(60,100)(14,14){5}{\line(1,1){10}}
\multiput(60,100)(14,-14){5}{\line(1,-1){10}}
\multiput(240,100)(14,14){5}{\line(1,1){10}}
\multiput(240,100)(14,-14){5}{\line(1,-1){10}}
\put(-10,166){\footnotesize{$\pi\kappa$}}
\put(-10,175){\footnotesize{Stokes line}}
\put(-10,30){\footnotesize{$-\pi\kappa$}}
\put(-10,20){\footnotesize{Stokes line}}
\put(90,166){\footnotesize{$\pi\kappa/2$}}
\put(85,26){\footnotesize{$-\pi\kappa/2$}}
\put(100,122){\footnotesize{$\pi\sigma$}}
\put(100,73){\footnotesize{$-\pi\sigma$}}
\put(-5,-10){$(a)$}
\put(170,166){\footnotesize{$\pi\kappa$}}
\put(170,175){\footnotesize{Stokes line}}
\put(170,30){\footnotesize{$-\pi\kappa$}}
\put(170,20){\footnotesize{Stokes line}}
\put(270,166){\footnotesize{$\pi\kappa/2$}}
\put(265,26){\footnotesize{$-\pi\kappa/2$}}
\put(175,-10){$(b)$}
\put(108,116){\circle*{4}}
\put(108,84){\circle*{4}}
\put(200,140){\circle*{4}}
\put(200,60){\circle*{4}}
\end{picture}
\caption{\small{The Stokes and anti-Stokes lines when $0<\kappa<1$ and the location of the arguments $P_\pm$ (indicated by heavy dots):
(a) $P_\pm=xe^{\pm\pi i\sigma}$ and (b) $P_\pm=xe^{\pm\pi i\kappa}$.}}
\end{figure}
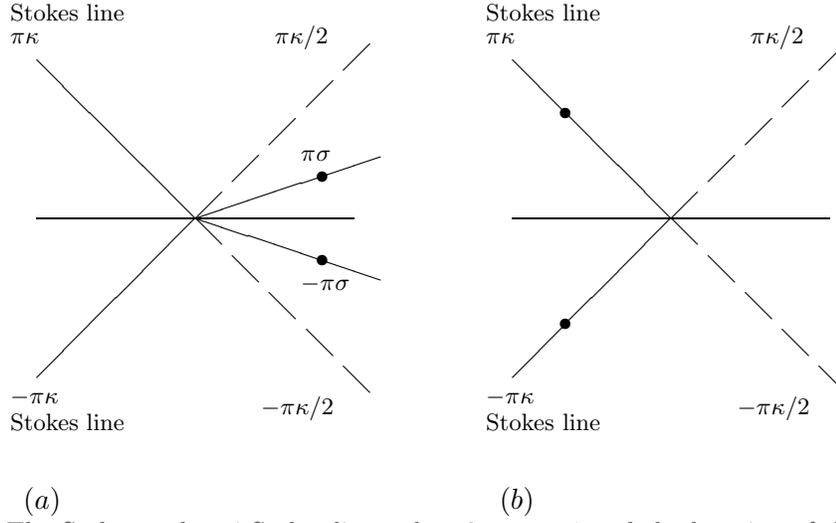

Let us denote the points $xe^{\pm\pi i\sigma}$ in the $\zeta$-plane that appear in the arguments of the associated function ${}_2\Psi_1(\zeta)$ in (\ref{e41}) by $P_\pm$. Then when
$0<\sigma<\f{1}{3}$, $P_\pm$ lie in the exponentially large sector $|\arg\,\zeta|<\fs\pi\kappa$ in Fig.~2(a) and consequently
${\hat E}_{2,1}(x)$ is exponentially large as $x\to+\infty$. When $\sigma=\f{1}{3}$, $P_\pm$ lie on the anti-Stokes lines $\arg\,\zeta=\pm\fs\pi\kappa$; on these rays $\cos \pi\sigma/\kappa=0$ and ${\hat E}_{2,1}(x)$ is oscillatory with an algebraically controlled amplitude. When $\f{1}{3}<\sigma<\fs$, $P_\pm$ lie in the exponentially small and algebraic sectors so that ${\hat E}_{2,1}(x)$ is exponentially small. When $\sigma=\fs$, $P_\pm$ lie on the Stokes lines $\arg\,\zeta=\pm\pi\kappa$, where the subdominant exponential expansions are in the process of switching off. Finally, when $\fs<\sigma<1$, $P_\pm$ are situated in the algebraic sectors, where the expansions are purely algebraic.

Then we obtain the following theorem.
\begin{theorem}$\!\!\!.$\ When $a=-\sigma$, with $0<\sigma<1$, we have the expansion of the generalised Struve function
\bee\label{e47}
(\fs z)^{-\nu-1} {\bf L}_\nu(z;-\sigma)\sim\left\{\begin{array}{ll} {\hat E}_{2,1}(x)+{\hat H}_{2,1}(x) & (0<\sigma<\fs) \\
\\
{\hat H}_{2,1}(x) & (\fs<\sigma<1)\end{array}\right.
\ee
as $z\to+\infty$, where ${\hat E}_{2,1}(x)$ and ${\hat H}_{2,1}(x)$ are defined in (\ref{e43}) and (\ref{e44}) with $x=z^2/4$, $X=\kappa(hx)^{1/\kappa}$. The quantities $\kappa$, $h$ and $A_0$ are defined in (\ref{e41a}) with $\vartheta=-\nu-\f{3}{2}$. The first few coefficients $c_j=c_j(-\sigma,\nu)$ are given in (\ref{e31a}).
\end{theorem}
When $\sigma=\fs$, the exponential expansion ${\hat E}_{2,1}(x)$ is maximally subdominant and is in the process of switching off. Thus, if we neglect this exponentially small contribution, the expansion of $(\fs z)^{-\nu-1} {\bf L}_\nu(z;-\fs)$ as $z\to+\infty$ is given by ${\hat H}_{2,1}(x)$.
\bigskip

\noindent 4.2 \ {\it The expansion for $z\to i\infty$}
\vspace{0.2cm}

\noindent
When $z=i|z|$, we have upon replacing $x$ by $xe^{\pi i}$ (so that $x=|z|^2/4$) and using the result ${}_2\Psi_1(\zeta e^{2\pi i})={}_2\Psi_1(\zeta)$
\bee\label{e48}
(\fs i|z|)^{-\nu-1} {\bf L}_\nu(i|z|;-\sigma)=\frac{i}{2\pi}\bl\{e^{\pi i\vartheta} {}_2\Psi_1(xe^{-\pi i\kappa})-e^{-\pi i\vartheta} {}_2\Psi_1(xe^{\pi i\kappa})\br\}.
\ee 
Then, for $x>0$, the associated functions ${}_2\Psi_1(xe^{\pm\pi i\kappa})$ have arguments situated on the Stokes lines, where the exponential contribution is maximally subdominant; see Fig.~2(b).

The algebraic contribution in (\ref{e48}) is
\begin{eqnarray*}
&&\frac{i}{2\pi}\bl\{e^{\pi i\vartheta} H_{2,1}(xe^{-\pi i\kappa} \cdot e^{\pi i})-e^{-\pi i\vartheta} H_{2,1}(xe^{\pi i\kappa} \cdot e^{-\pi i})\br\}\\
&=&-\frac{1}{\pi}\sum_{k=0}^\infty (-)^k G_k x^{-k-1} \sin \pi(\vartheta-\sigma(k+1))-\frac{1}{\pi\sigma}\sum_{k=0}^\infty (-)^k G_k' x^{-k_s} \sin \pi(\vartheta-\sigma k_s)\\
&=&\frac{1}{\pi}\sum_{k=0}^\infty \frac{\g(k+\fs) x^{-k-1}}{\g(\nu+\f{3}{2}+\sigma(k+1))},
\end{eqnarray*}
the second sum vanishing since $\sin \pi(\vartheta-\sigma k_s)=\sin \pi k \equiv 0$ for $k=0, 1, 2, \ldots\ $.

Then we have
\begin{theorem}$\!\!\!.$\ When $a=-\sigma$, with $0<\sigma<1$, we have the expansion of the generalised Struve function
\bee\label{e49}
(\fs i|z|)^{-\nu-1} {\bf L}_\nu(i|z|;-\sigma)\sim\frac{1}{\pi}\sum_{k=0}^\infty \frac{\g(k+\fs) (\fs|z|)^{-2k-2}}{\g(\nu+\f{3}{2}+\sigma(k+1))}
\ee
as $|z|\to\infty$. Here we have neglected the maximally subdominant exponentially small contribution to the expansion.
\end{theorem}
\vspace{0.6cm}

\begin{center}
{\bf 5. \ Numerical results}
\end{center}
\setcounter{section}{5}
\setcounter{equation}{0}
\renewcommand{\theequation}{\arabic{section}.\arabic{equation}}
Some numerical examples are presented to illustrate the expansions developed in Sections 3 and 4. We first consider the expansions valid when $a>0$. In Table 1 we show the normalised coefficients $c_j$ corresponding to $a=\fs$ and $\nu=\f{1}{4}$ computed by means of the algorithm described in \cite[Appendix A]{P17a}. Table 2 shows values\footnote{In the tables we write $x(y)$ to represent $x\times 10^y$.} of the generalised Struve function when $a>0$. The top half of the table gives values of $${\cal L}_\nu(z;a)\equiv (\fs z)^{_\nu-1} {\bf L}_\nu(z;a)$$ for different $z>0$ when $a=\fs$, $\nu=\f{1}{4}$
compared with the exponential expansion $E_{1,2}(\zeta)$ in (\ref{e33a}). In the computation of $E_{1.2}(\zeta)$ we have employed the truncation index $j=10$. The lower half of the table gives values of ${\cal L}_\nu(i|z|;a)$ with the dominant, optimally truncated\footnote{That is, truncated at, or near, the term of smallest magnitude.} algebraic expansion $H_{1,2}(x)$ (with $x=|z|^2/4$) subtracted off. This value is compared with the exponential expansion ${\hat E}_1(X)$ defined in (\ref{e33}).
\begin{table}[th]
\caption{\footnotesize{The normalised coefficients $c_j=A_j/A_0$ for ${}_1\Psi_2(\zeta)$ in (\ref{e31}) when $a=\fs$ and $\nu=\f{1}{4}$. }}
\begin{center}
\begin{tabular}{c|l||c|l}
\mcol{1}{c|}{$j$} & \mcol{1}{c||}{$c_j$} & \mcol{1}{c|}{$j$} & \mcol{1}{c}{$c_j$}\\
[.05cm]\hline
&&& \\[-0.2cm]
1 & $-\f{5}{12}$ & 2 & $-\f{35}{288}$\\
&&& \\[-0.2cm]
3 & $-\f{665}{10368}$& 4 & $+\f{9625}{497664}$\\
&&& \\[-0.2cm]
5 & $+\f{1856855}{5971968}$ & 6 & $+\f{606631025}{429981696}$\\
&&& \\[-0.2cm]
7 & $+\f{27773871125}{5159780352}$ & 8 & $+\f{8996211899675}{495338913792}$\\
&&& \\[-0.2cm]
9 & $+\f{2459153764892825}{53496602689536}$ & 10 & $-\f{22173972436540925}{1283918464548864}$\\
[.15cm]\hline
\end{tabular}
\end{center}
\end{table}

\begin{table}[th]
\caption{\footnotesize{Values of ${\cal L}_\nu(z;a)\equiv (\fs z)^{-\nu-1}{\bf L}_\nu(z;a)$ compared with the asymptotic expansions in (\ref{e33a}) and (\ref{e33}) when $a=\fs$ and $\nu=\f{1}{4}$. The exponential expansions have truncation index $j\leq 10$ and the algebraic expansion $H_{1,2}(x)$ (with $x=|z|^2/4$) has been optimally truncated.}}
\begin{center}
\begin{tabular}{c|ll}
\mcol{1}{c|}{$z$} & \mcol{1}{c}{${\cal L}_\nu(z;a)$} & \mcol{1}{c}{$E_{1,2}(\zeta)$}\\
[.05cm]\hline
&& \\[-0.2cm]
5  & $+3.452097942\,(1)$ & $+3.461544352\,(1)$\\
10 & $+1.226039040\,(5)$ & $+1.226039286\,(5)$\\
12 & $+6.877617187\,(6)$ & $+6.877617204\,(6)$\\
15 & $+5.182624938\,(9)$ & $+5.182624938\,(9)$\\ 
[.15cm]\hline
&& \\[-0.3cm]
\mcol{1}{c|}{$z$} & \mcol{1}{c}{${\cal L}_\nu(z;a)-H_{1,2}(x)$} & \mcol{1}{c}{${\hat E}_1(X)$}\\
[.05cm]\hline
&& \\[-0.2cm]
10$i$ & $-4.572292174\,(-6)$  & $-4.57195324?\,(-6)$\\
15$i$ & $+4.021530098\,(-10)$ & $+4.021543491\,(-10)$\\
20$i$ & $+6.827666326\,(-12)$ & $+6.827666325\,(-12)$\\
25$i$ & $+3.515426867\,(-15)$ & $+3.515426867\,(-15)$\\
[.15cm]\hline
\end{tabular}
\end{center}
\end{table}

\begin{table}[th]
\caption{\footnotesize{Values of ${\cal L}_\nu(z;-\sigma)\equiv (\fs z)^{-\nu-1}{\bf L}_\nu(z;-\sigma)$ compared with the asymptotic expansions in (\ref{e47}) when $z>0$ and $\nu=\f{1}{3}$. The exponential expansion has been optimally truncated.}}
\begin{center}
\begin{tabular}{r|ll||r|ll}
\mcol{1}{c|}{} & \mcol{2}{c||}{$\sigma=1/5$} & \mcol{1}{c|}{} & \mcol{2}{c}{$\sigma=1/3$}\\
\mcol{1}{c|}{$z$} & \mcol{1}{c}{${\cal L}_\nu(z;-\sigma)$} & \mcol{1}{c||}{Asymptotic} & \mcol{1}{c|}{$z$} & \mcol{1}{c}{${\cal L}_\nu(z;-\sigma)$} & \mcol{1}{c}{Asymptotic}\\
[.05cm]\hline
&&&& \\[-0.25cm]
8  & $+1.371278215\,(4)$ & $\!\!\!\! +1.371994397\,(4)$ & 5 &  $+4.994707877\,(1)$ & $\!\!\!\!+4.992261627\,(1)$\\
10 & $-1.628234940\,(7)$ & $\!\!\!\!-1.628235076\,(7)$  & 8 &  $+5.127188845\,(2)$ & $\!\!\!\!+5.127188845\,(2)$\\
15 & $-2.287676991\,(22)$ & $\!\!\!\!-2.287676991\,(22)$& 10&  $+1.563077837\,(3)$ & $\!\!\!\!+1.563077837\,(3)$\\
\hline
&&&& \\[-0.3cm]
\mcol{1}{c|}{} & \mcol{2}{c||}{$\sigma=1/2$} & \mcol{1}{c|}{} & \mcol{2}{c}{$\sigma=3/5$}\\
\mcol{1}{c|}{$z$} & \mcol{1}{c}{${\cal L}_\nu(z;-\sigma)$} & \mcol{1}{c||}{Asymptotic} & \mcol{1}{c|}{$z$} & \mcol{1}{c}{${\cal L}_\nu(z;-\sigma)$} & \mcol{1}{c}{Asymptotic}\\
[.05cm]\hline
&&&& \\[-0.25cm]
3  & $+4.705453951\,(0)$ & $\!\!\!\!+4.719691159\,(0)$ & 3 &  $+4.075339511\,(0)$ & $\!\!\!\!+4.074935642\,(0)$\\
5  & $+1.918197617\,(1)$ & $\!\!\!\!+1.918197638\,(1)$ & 4 &  $+7.439302510\,(0)$ & $\!\!\!\!+7.439299037\,(0)$\\
8  & $+8.747082153\,(1)$ & $\!\!\!\!+8.747082153\,(1)$ & 5 &  $+1.276299496\,(1)$ & $\!\!\!\!+1.276299496\,(1)$\\
\hline
\end{tabular}
\end{center}
\end{table}
\begin{table}[th]
\caption{\footnotesize{Values of ${\cal L}_\nu(z;-\sigma)\equiv (\fs z)^{-\nu-1}{\bf L}_\nu(z;-\sigma)$ compared with the asymptotic expansion (\ref{e49}) when $\arg\,z=\fs\pi$ and $\nu=\f{4}{3}$. The algebraic expansion (\ref{e49}) has been optimally truncated.}}
\begin{center}
\begin{tabular}{r|ll||r|ll}
\mcol{1}{c|}{} & \mcol{2}{c||}{$\sigma=1/4$} & \mcol{1}{c|}{} & \mcol{2}{c}{$\sigma=1/3$}\\
\mcol{1}{c|}{$z$} & \mcol{1}{c}{${\cal L}_\nu(z;-\sigma)$} & \mcol{1}{c||}{Asymptotic} & \mcol{1}{c|}{$z$} & \mcol{1}{c}{${\cal L}_\nu(z;-\sigma)$} & \mcol{1}{c}{Asymptotic}\\
[.05cm]\hline
&&&& \\[-0.25cm]
6$i$ &  $3.044656205\,(-2)$ & $3.044653596\,(-2)$ & 6$i$ &  $2.792844201\,(-2)$ & $2.792844405\,(-2)$\\
8$i$ &  $1.673275565\,(-2)$ & $1.673275565\,(-2)$ & 8$i$ &  $1.539185802\,(-2)$ & $1.539185802\,(-2)$\\
\hline
&&&& \\[-0.3cm]
\mcol{1}{c|}{} & \mcol{2}{c||}{$\sigma=1/2$} & \mcol{1}{c|}{} & \mcol{2}{c}{$\sigma=3/4$}\\
\mcol{1}{c|}{$z$} & \mcol{1}{c}{${\cal L}_\nu(z;-\sigma)$} & \mcol{1}{c||}{Asymptotic} & \mcol{1}{c|}{$z$} & \mcol{1}{c}{${\cal L}_\nu(z;-\sigma)$} & \mcol{1}{c}{Asymptotic}\\
[.05cm]\hline
&&&& \\[-0.25cm]
4$i$  & $5.552864403\,(-2)$ & $5.553062223\,(-2)$ & 3$i$ &  $7.704243224\,(-2)$ & $7.704358006\,(-2)$\\
5$i$  & $3.420993477\,(-2)$ & $3.420993479\,(-2)$ & 4$i$ &  $4.087728092\,(-2)$ & $4.087728092\,(-2)$\\
\hline
\end{tabular}
\end{center}
\end{table}

We now consider the expansions valid when $-1<a<0$.
Table 3 shows values of ${\cal L}_\nu(z;-\sigma)$ for different $z>0$ and $\sigma$ in the range $0<\sigma<1$ compared with the asymptotic expansion in (\ref{e47}). For each value of $\sigma$ the normalised coefficients $c_j$ were computed. The expansions have been optimally truncated and we have used the result ${\cal L}_\nu(z;-\sigma)\sim H_{2,1}(x)$ as $z\to+\infty$ (with $x=z^2/4$) when $\sigma=\fs$. Finally, in Table 4 we show values of ${\cal L}_\nu(i|z|;-\sigma)$ for different $|z|$ and $\sigma$ in the range $0<\sigma<1$ compared with the optimally truncated algebraic expansion in (\ref{e49}).

It is seen in all cases that very good agreement is achieved when $z$ is sufficiently large.

\vspace{0.6cm}

\begin{center}
{\bf Appendix: \ The coefficients $c_j\equiv c_j(-\sigma,\nu)$ in the expansion (\ref{e43})}
\end{center}
\setcounter{section}{1}
\setcounter{equation}{0}
\renewcommand{\theequation}{\Alph{section}.\arabic{equation}}
The inverse factorial expansion (\ref{e22a}) involving the coefficients $c_j\equiv c_j(\alpha,\nu)$ for the function ${}_1\Psi_2$ in (\ref{e31}) can be written in the form (see \cite[Appendix A]{P17a})
\[R(\alpha)\equiv\frac{\g(\kappa s+\vartheta')}{\g(\alpha s+\nu+\f{3}{2})\g(s+\f{3}{2})}=\kappa A_0 (h\kappa^\kappa)^s\bl\{\sum_{j=0}^{M-1} \frac{c_j(\alpha,\nu)}{(\kappa s+\vartheta')_j}+\frac{O(1)}{(\kappa  s+\vartheta')_M}\br\}\]
as $|s|\to\infty$ in $|\arg\,s|<\pi$, where $\kappa=1+\alpha$, $h=\alpha^{-\alpha}$, $\vartheta'=-\nu-\f{3}{2}$ and $A_0=(\kappa/\alpha)^{\nu+1}/\sqrt{2\pi}$. Then, when $\alpha=-\sigma$, $0<\sigma<1$, we have
\[R(-\sigma)=\frac{\kappa}{\sqrt{2\pi}} \bl(\!\frac{\kappa}{-\sigma}\!\br)^{\nu+1} ((-\sigma)^\sigma \kappa^\kappa)^s
\bl\{\sum_{j=0}^{M-1} \frac{c_j(-\sigma,\nu)}{(\kappa s+\vartheta')_j}+\frac{O(1)}{(\kappa s+\vartheta')_M}\br\}\]
\[=\frac{\kappa}{\sqrt{2\pi}}\, \bl(\frac{\kappa}{\sigma}\br)^{\nu+1} (\sigma^\sigma \kappa^\kappa)^s
e^{\pi i(\sigma s-\nu-1)}
\bl\{\sum_{j=0}^{M-1} \frac{c_j(-\sigma,\nu)}{(\kappa s+\vartheta')_j}+\frac{O(1)}{(\kappa a+\vartheta')_M}\br\}.\]

The inverse factorial expansion (\ref{e22a}) for the function ${}_2\Psi_1$ in (\ref{e41b}) becomes (with $\kappa$, $h$ and $A_0$ defined in (\ref{e41a}))
\[S\equiv \frac{\g(\kappa s+\vartheta') \g(\sigma s-\nu-\fs)}{\g(s+\f{3}{2})}\]
\bee\label{a1}
=\kappa \sqrt{2\pi} \bl(\frac{\kappa}{\sigma}\br)^{\nu+1} (\sigma^\sigma \kappa^\kappa)^s
\br\{\sum_{j=0}^{M-1} \frac{d_j}{(\kappa s+\vartheta')_j}+\frac{O(1)}{(\kappa s+\vartheta')_M}\br\},
\ee
where $d_j$ are coefficients to be determined. But
\[S=\frac{\pi R(-\sigma)}{\sin \pi(\sigma s-\nu-\fs)}=\frac{2\pi R(-\sigma)}{e^{\pi i(\sigma s-\nu-1)} \Lambda(s)}\]
\bee\label{a2}=\frac{\kappa \sqrt{2\pi}}{\Lambda(s)}\,\bl(\frac{\kappa}{\sigma}\br)^{\nu+1} (\sigma^\sigma \kappa^\kappa)^s
\bl\{\sum_{j=0}^{M-1} \frac{c_j(-\sigma,\nu)}{(\kappa s+\vartheta')_j}+\frac{O(1)}{(\kappa a+\vartheta')_M}\br\},
\ee
where $\Lambda(s)=1-e^{-2\pi i(\sigma s-\nu-\frac{1}{2})}$.

Letting $|s|\to\infty$ in the sector $-\pi<\arg\,s<0$, so that $\Lambda(s)\to 1$, we have upon comparison of the terms in
(\ref{a1}) and (\ref{a2}) that $d_j=c_j(-\sigma,\nu)$ \ ($1\leq j\leq M-1$).

\end{document}